\documentclass[12pt]{extarticle}
\usepackage{amsmath, amsthm, amssymb, hyperref, color}
\usepackage{graphicx}
\usepackage{caption}
\usepackage{subcaption}
\usepackage{mathtools}
\usepackage{enumerate}
\usepackage{stackengine}
\usepackage[all]{xypic}
\usepackage{verbatim}
\usepackage{upgreek}
\usepackage{chemfig, chemnum}
\usepackage{tikz}
\usepackage[lined,commentsnumbered,ruled]{algorithm2e}
\usepackage{caption}
\usepackage{subcaption}
\usepackage{float}
\oddsidemargin \evensidemargin
\hypersetup{
    colorlinks=false,
    linkcolor=blue,
    filecolor=blue,      
    urlcolor=blue,
}
\theoremstyle{theorem}
\newtheorem{theorem}{Theorem}[section]

\newtheorem{lemma}[theorem]{Lemma}

\theoremstyle{definition}

\newtheorem{remark}[theorem]{Remark}

\newtheorem{example}[theorem]{Example}

\newcommand{\CC}{\mathbb{C} }

\title{\textbf{Numerical reconstruction of curves from their Jacobians}}

\author{{Daniele Agostini, T\"urk\"u \"Ozl\"um \c{C}elik, Demir Eken}}

\date{}

\begin{document}

\maketitle

\begin{abstract}
\noindent 
We approach the Torelli problem of recostructing a curve from its Jacobian from a computational point of view. Following Dubrovin, we design a machinery to solve this problem effectively, which builds on methods in numerical algebraic geometry. We verify this methods via numerical experiments with curves up to genus 7.
\end{abstract}

\section{Introduction}
\label{introduction}

The Torelli theorem is a classical and foundational result in algebraic geometry, stating that a Riemann surface, or smooth algebraic curve, $C$ is uniquely determined by its Jacobian variety $J(C)$. More concretely, the theorem says that a Riemann surface of genus $g$ can be recovered from one Riemann matrix $\tau$ that represents its Jacobian, together with the theta function
\begin{equation}\label{eq:theta}
\theta\colon \mathbb{C}^g \times \mathbb{H}_g \longrightarrow \mathbb{C}, 
\qquad \theta({\bf z},\tau) := \sum_{n\in \mathbb{Z}^g} \exp\left( \pi i n^t\tau n + 2\pi i n^t{\bf z}  \right)
\end{equation}
where $\mathbb{H}_g$ is the Siegel upper-half space of $g\times g$ symmetric complex matrices with positive definite imaginary part. There are various proofs of Torelli's theorem, which can be even made concrete in computational terms. Most proofs rely on the geometry of the theta divisor. This is the locus inside the Jacobian variety $J(C) = \mathbb{C}^g / \mathbb{Z}^g + \tau\mathbb{Z}^g$ which is cut out by the theta function:
\[ \Theta = \{ {\bf z}\in J(C) \,|\, \theta({\bf z},\tau) = 0  \}. \]
For example, suppose that the Riemann surface $C$ is not hyperelliptic, so that we can identify $C$ with a canonical model $C\subseteq \mathbb{P}^{g-1}$.  Then for any singular point ${\bf z}\in \Theta_{\rm sing}$ of the theta divisor, the corresponding Hessian matrix
$(\frac{\partial^2 \theta}{\partial z_i \partial z_j}({\bf z},\tau))$
defines a quadric in the projective space space $\mathbb{P}^{g-1}$. By a result of Green \cite{Green}, such quadrics span the space of quadrics in the ideal of the curve. Hence, if the curve is not trigonal, or a smooth plane quintic, these quadrics generate the whole canonical ideal. This result has been extended by Kempf and Schreyer, which gave a way to recover the curve from a single singular point \cite{KS}. In particular, this gives a powerful effective reconstruction of the curve, provided that we are able to solve the system
\begin{equation}\label{eq:singularpointstheta} 
\theta({\bf z},\tau) = \frac{\partial \theta}{\partial z_1 }(z,\tau) = \dots = \frac{\partial \theta}{\partial z_g}({\bf z},\tau) = 0. 
\end{equation}
This has been implemented numerically for curves of genus $4$ in \cite{ChuKumStu}, but it is a rather hard task in general, since the theta function is inherently transcendental. Moreover, this problem is also quite sensitive to the precision of the data: for example, if we move $\tau$ a bit, the corresponding theta divisor will not have singular points. There are various other proofs of the Torelli theorem, but many involve solving system of equations such as \eqref{eq:singularpointstheta}. 

Hence, we look for different, more algebraic methods. Such a strategy was proposed by Dubrovin \cite{Dub81}, building on Krichever's work \cite{Kri} on algebraic curves and the Kadomtsev-Petviashvili (KP) equation:
\begin{equation}
\label{eq:KP1}
\frac{\partial}{\partial x}\left( 4u_t - 6uu_x -u_{xxx} \right)\, = \, 3u_{yy}.
\end{equation}
More precisely, for each Riemann surface $C$ of genus $g$ there exists a threefold $\mathcal{D}_C$ in a weighted projective space $\mathbb{WP}^{3g-1}$ parametrizing triples $(U,V,W)$ such that the function
\begin{equation}
\label{eq:u_tau}
u(x,y,t) \,\,=\,\, 2\frac{\partial^2}{\partial x^2} \log \uptau(x,y,t)+c, \qquad \uptau(x,y,z) := \theta(Ux+Vy+Wz+D,\tau)
\end{equation}
is a solution to the KP equation \eqref{eq:KP1} for any $D\in \mathbb{C}^g$ and some $c\in \mathbb{C}$. This threefold was called the \emph{Dubrovin threefold} in \cite{AgoCelStu} and it was studied there from a computational point of view. The important properties of this object for our point of view are two: first, $\mathcal{D}_C$ is cut out by some explicit equations whose coefficient are derivatives of theta functions (with characteristic) evaluated at zero. These can be computed explicitly with a software for the evaluation of the theta functions, such as \texttt{Theta.jl} in Julia \cite{julia}. Second, the projection of $\mathcal{D}_C$ onto the projective space $\mathbb{P}^{g-1}_U$ of the coordinates $u_1,\dots,u_g$ consists exactly of the canonical model for the curve $C\subseteq \mathbb{P}^{g-1}_U$.
Hence, equations for the canonical model of $C$ can be obtained by eliminating the variables $V,W$ from the equations of the Dubrovin threefold $\mathcal{D}_C$, a purely algebraic process. In conclusion, this allows to recover the curve from the Riemann matrix $\tau$ without having to solve a transcendental system such as \eqref{eq:singularpointstheta}. 

In this note, we explain how to implement this strategy effectively, using the methods of numerical algebraic geometry. In Section \ref{sec:Dubrovin}, we explain the background behind the Dubrovin threefold and we state the key Lemma \ref{lemma:quarticelimination}, which explains how to recover equations for the curve. In particular, this allows to recover quartics equations, but we discuss also the case of quadrics and cubics. Furthermore, we also comment on applications of these methods to the classical Schottky problem. In Section \ref{dubrovin} we state the algorithm and we analyze its complexity. Moreover, even if our focus is on methods that avoid finding singular points, this can be very useful when we are able to solve the transcendental system \eqref{eq:singularpointstheta}, and we comment on this in Section \ref{SingularPointsTheta}. We conclude by presenting numerical experiments with curves from genera from  3 to 7, which we carried out with the packages \texttt{RiemannSurfaces} in Sage and \texttt{Theta.jl} and \texttt{Homotopycontinuation.jl} in Julia.

\section{The Dubrovin Threefold}\label{sec:Dubrovin}

We start by recalling some background on the Dubrovin threefold, follwing \cite{AgoCelStu,Dub81}. Let $C$ be a smooth projective algebraic curve, or compact Riemann surface, of genus $g$. We fix a symplectic basis $\,a_1,b_1,\dots,a_g,b_g\,$ for the
first homology group $H_1(C,\mathbb{Z})$ and we choose a \emph{normalized basis} $\omega_1,\omega_2,\dots,\omega_g$ of holomorphic differentials, meaning that 
\begin{equation}\label{eq:adaptedbasis} 
\int_{a_i} \omega_j\,\,=\,\, \delta_{ij} ,
\end{equation}
The corresponding Riemann matrix $\tau \in \mathbb{H}_g$ is defined as
\begin{equation}\label{eq:riemannmatrix}
\tau = \left( \int_{b_i}\omega_j \right)_{1\leq i,j \leq g}.
\end{equation}
One can see that the function \eqref{eq:u_tau} is a solution to the KP equation \eqref{eq:KP1} for all values of $D\in \mathbb{C}^g$ and a certain value of $c\in \mathbb{C}$ if and only if there exists a $d\in \mathbb{C}$ such that the following quartic PDE, known as the Hirota bilinear equation, is satisfied:
\begin{equation}\label{eq:hirota}
\!\! ( \uptau_{xxxx}\uptau -4\uptau_{xxx}\uptau_x + 3\uptau_{xx}^2)
\,+\,4 (\uptau_{x}\uptau_t \,-\,  \uptau \uptau_{xt}) \,+\, 6 c\, ( \uptau_{xx}\uptau \, -\, \uptau_x^2)
\,+\, 3 (\uptau \uptau_{yy}\,-\, \uptau_y^2)\,+\,8d\,\uptau^2 
\,\,= \,\,0.
\end{equation}
We now introduce a weighted projective space $\mathbb{WP}^{3g+1}$ with variables $(U,V,W,c,d)$ where the $U=(u_1,\dots,u_g)$ have degree $1$, the $V=(v_1,\dots,v_g)$ have degree $2$, the $W=(w_1,\dots,w_g)$ have degree $3$ and finally $c,d$ have degree $2$ and $4$ respectively. The \emph{big Dubrovin threefold} $\mathcal{D}_{C}^{\rm big}$ parametrizes all elements $(U,V,W,c,d)$, with $U\ne \mathbf{0}$, such that $\uptau(x,y,z)$ in \eqref{eq:u_tau} is a solution to the Hirota bilinear equation \eqref{eq:hirota} for all $D\in \mathbb{C}^g$. The projection of this variety to the space $\mathbb{WP}^{3g-1}$ of the $(U,V,W)$ is called simply the \emph{Dubrovin threefold} $\mathcal{D}_{C}$.

Equations for $\mathcal{D}^{\rm big}_C$ can be obtained directly from \eqref{eq:hirota} as follows.  Given any ${\bf z}$ in $\CC^g$, we write the Riemann theta function as $\theta({\bf z}) = \theta({\bf z} , \tau)$. Then we consider the differential operator $\partial_U := u_1 \frac{\partial}{\partial z_1} + \cdots + u_g \frac{\partial}{\partial z_g}$, and the analogous operators $\partial_V, \partial_W$. For any fixed vector $\mathbf{z}\in\mathbb{C}^g$, the \emph{Hirota quartic} $H_{\bf z}$ is defined as:
\begin{multline}\label{eq:hirotaquartics}
\left(\partial_U^4\theta({\bf z}) \cdot \theta({\bf z})-4\partial^3_U\theta({\bf z})\cdot \partial_U\theta({\bf z})+3\{\partial^2_U\theta({\bf z})\}^2\right)+4\cdot \left(\partial_U\theta({\bf z})\cdot \partial_W\theta({\bf z})-\theta({\bf z}) \cdot\partial_U\partial_W \theta({\bf z})\right)\\+6c\cdot \left( \partial^2_U\theta({\bf z})\cdot \theta({\bf z}) - \{\partial_U\theta({\bf z})\}^2 \right) + 3\cdot \left(\theta({\bf z})\cdot \partial^2_V\theta({\bf z})-\{\partial_V\theta({\bf z})\}^2 \right) +8d\cdot \theta({\bf z})^2.
\end{multline}
This is exactly the expression obtained by combining \eqref{eq:hirota} and \eqref{eq:KP1}, hence the big Dubrovin threefold $\mathcal{D}^{\rm big}_C$ is cut out by  the
Hirota quartics $H_{\bf z}$, as $\mathbf{z}$ runs over all vectors in $ \mathbb{C}^g$, see \cite[Proposition 4.2]{AgoCelStu}. The coefficients $H_{\bf z}(U,V,W,c,d)$ are the values of 
the theta function $\theta$ and its partial derivatives of certain order at ${\bf z}$ and they can be computed using numerical software for evaluating theta functions and their derivatives. We use the {\tt Julia} package that is introduced in \cite{julia}.

This yields an infinite number of equations for the Dubrovin threefold. A finite number can be obtained via the theta functions with characteristic $\varepsilon,\delta \in \{ 0,1\}^g$: 
 \begin{equation}
\label{eq:thetaFunctionChar}
\theta\begin{bmatrix} \varepsilon \\ \delta \end{bmatrix}({\bf z}\, |\, \tau)\,\,\, = \,\,\,
\sum_{n \in \mathbb{Z}^g} {\rm exp} \left( \pi i \left(n+\frac{\varepsilon}{2}\right)^T \tau \left(n+\frac{\varepsilon}{2}\right) + \left(n + \frac{\varepsilon}{2}\right)^T  \left(z + \frac{\delta}{2}\right) \right).
\end{equation}
This coincides with the Riemann theta function \eqref{eq:theta} for $\varepsilon=\delta=0$ and in general it differs from it by an exponential factor.  We consider the following function
\begin{equation}
\label{eq:doubled}
\hat{\theta}[\varepsilon]({{\bf z}})\,\, :=\,\, {\theta}\begin{bmatrix} \varepsilon \\ 0 \end{bmatrix}({\bf z}\, |\, 2\tau) .
\end{equation}
For fixed $\tau$, these complex numbers $\,\hat{\theta}[\varepsilon](0)\,$ at ${\bf z}=0$ are called \emph{theta constants}.
We use the term theta constant also for evaluations at ${\bf z}={\bf 0}$ of derivatives of (\ref{eq:doubled}).
With these conventions, we define the \emph{Dubrovin quartic} in $(U,V,W,c,d)$ associated to the half-characteristic $\varepsilon$ as:
\begin{equation}\label{eq:dubrovinquartic}
F[\varepsilon](U,V,W,c,d) \,\,:=\,\, \partial^4_U \hat{\theta}[\varepsilon]({\bf 0})-\partial_U\partial_{W}\hat{\theta}[\varepsilon]({\bf 0})+\frac{3}{2}c\cdot \partial^2_U\hat{\theta}[\varepsilon]({\bf 0})+\frac{3}{4}\partial^2_V\hat{\theta}[\varepsilon]({\bf 0})+d\hat{\theta}[\varepsilon]({\bf 0}).
\end{equation}
%We call $F[\varepsilon](U,V,W,c,d)$  the {\em Dubrovin quartic} associated with the half-characteristic $\varepsilon$.
The Dubrovin and the Hirota quartics span  the same vector subspace of $\,\CC[U,V,W,c,d\,]_4$ as shown in \cite[Proposition 4.3]{AgoCelStu}, hence they also provide defining equations for $\mathcal{D}_C^{\rm big}$.

We come to the crucial point: the projection of the big Dubrovin threefold onto the projective space $\mathbb{P}^{g-1} = \mathbb{P}^{g-1}_{U}$ coincides exactly with the canonical model of the curve $C$ induced by the basis of holomorphic differentials of \eqref{eq:adaptedbasis}:
\begin{equation}\label{eq:canonicalmodel}
C \longrightarrow \mathbb{P}^{g-1}_U, \qquad p \mapsto \left[ \omega_1(p),\omega_2(p),\dots,\omega_g(p) \right]
\end{equation}
in particular, if the curve $C$ is not hyperelliptic, the canonical model is isomorphic to the curve $C$ itself. In algebraic terms, this means that the canonical model of \eqref{eq:canonicalmodel} can be recovered by eliminating the variables $V,W,c,d$ from the equations of the big Dubrovin threefold.  This is reduced to a problem of linear algebra as follows: for any half-characteristic $\varepsilon \in \{0,1 \}^g$ write $Q[\varepsilon]$ for the Hessian matrix of the function $\hat{\theta}({\bf z})$, then, combining \cite[Lemma 4.6]{AgoCelStu} and \cite[Proposition 4.7]{AgoCelStu} we have:

\begin{lemma}\label{lemma:quarticelimination}
	Let us denote by $V_{\tau} \subseteq \mathbb{C}[u_1,\dots,u_g]$ the vector space of linear combinations
	\begin{equation}\label{eq:quarticelimination} 
	\sum_{\varepsilon \in \{0,1 \}^g} \lambda_{\varepsilon} \cdot \partial^4_U\hat{\theta}[\varepsilon], 
	\end{equation}
	where the $2^g$ complex scalars $\lambda_{\varepsilon}$ satisfy the linear equations
	\begin{equation}\label{eq:linearcondition}
	\sum_{\varepsilon} \lambda_{\varepsilon} \cdot Q[\varepsilon] \,=\, 0 \qquad {\rm and}  \qquad
	\sum_{\varepsilon} \lambda_{\varepsilon} \cdot \hat{\theta}[\varepsilon] \,=\, 0.
	\end{equation}
	Then a linear combination of the Dubrovin quartics is independent of $c,d$ if and only if it belongs to $V_{\tau}$. Furthermore, $V_{\tau}$ has dimension $\,2^g-\frac{g(g+1)}{2}-1$  and the corresponding quartics \eqref{eq:quarticelimination} cut out the canonical model \eqref{eq:canonicalmodel} of the curve $C$.
\end{lemma}

Hence, if the curve $C$ is not hyperelliptic, this Lemma gives a way to recover the curve from the Riemann matrix $\tau$, which depends only on the evaluation of the theta function and its derivatives.

\subsection{Recovering quadrics and cubics}\label{sec:quadrics}

Lemma \ref{lemma:quarticelimination} allows us to recover a linear space of quartics that cut out a canonical model of $C$. However, it is also possible to recover quadric equations.  We start with the following basic observation: if in the space $V_{\tau}$ we can find a quartic of the form $Q(U)^2$, then $Q$ is a quadric containing the curve $C$. We can actually find such special quartics inside $V_{\tau}$: indeed, suppose that $\mathbf{z}_0 \in \mathbb{C}^g$ is a singular point of the theta divisor $\Theta$. Then the Hirota quartic $H_{\bf z_0}$ becomes:
\begin{equation}\label{eq:hirotasingular}
 H_{\bf z_0} = 3\left( \partial^2_U \theta({\bf z_0})\right)^2
\end{equation}
and since this is independent of $c,d$, Lemma \ref{lemma:quarticelimination} tells us that $(\partial^2_U\theta({\bf z_0}))^2 \in V_{\tau}$. Furthermore, we know by Green's result mentioned in the introduction, that the quadrics $\partial^2_U\theta({\bf z_0})$ appearing in \eqref{eq:hirotasingular} span the whole vector space of quadrics in the ideal of the canonical curve. Hence, if $C$ is not hyperelliptic, trigonal, or a smooth plane quintic, such quadrics generate the canonical ideal of the curve. We again point out that, at least in principle, such quadrics can be computed by algebraic and not transcendental methods. Indeed, this corresponds to intersecting the space $V_{\tau}$ with the subvariety in $\mathbb{C}[U]_4$ given by quartics of the form $Q^2$, so it amounts to solving a polynomial system of equations in the space $V_{\tau}$. 

We discuss briefly also the case of cubics, which can appear if the curve is trigonal or a smooth plane quintic. In general, if $\mathbf{z}_0 \in \Theta$ is a singular point in the theta divisor, the cubic equation $\partial_U^3\theta({\bf z}_0)$ belongs to the canonical ideal of the curve \cite{KS}. If we apply the operator $\partial_U$ to the Hirota quartic $H_{\bf z}$ and we evaluate it at $\mathbf{z} = \mathbf{z}_0$, we obtain the quintic equation
%\begin{multline} \label{eq:delUhirotaQuart}
%{\partial_U H_{\bf z}} = \partial_U^4\theta({\bf z}) \cdot \theta({\bf z}) + \partial_U^5\theta({\bf z}) \cdot \partial_U\theta({\bf z}) - 4 (\partial_U^3\theta({\bf z}) \cdot \partial_U^2\theta({\bf z}) + \partial_U^4\theta({\bf z}) \cdot \partial_U\theta({\bf z}) ) + 6 \partial_U^2\theta({\bf z}) \cdot \partial_U^3\theta({\bf z}) \\+ 4 \cdot (\partial_U^2\theta({\bf z}) \cdot \partial_W\theta({\bf z}) + \partial_U\theta({\bf z}) \cdot \partial_U - (\partial_U\theta({\bf z}) \partial_W\theta({\bf z}) \cdot \partial_U \partial_W\theta({\bf z}) + \theta({\bf z}) \cdot \partial_U^2\partial_W\theta({\bf z}))) \\+ 6c \cdot (\partial_U^3\theta({\bf z}) \cdot \theta({\bf z})+ \partial_U^2\theta({\bf z}) \cdot \partial_U\theta({\bf z}) - 2 \partial_U\theta({\bf z}) \cdot \partial_U^2\theta({\bf z})) \\+ 3 \cdot (\partial_U\theta({\bf z}) \cdot \partial_V^2\theta({\bf z}) + \theta({\bf z}) \cdot \partial_U\partial_V^2\theta({\bf z}) - 2 \partial_V\theta({\bf z}) \cdot \partial_U\partial_V\theta({\bf z})) + 16d \cdot \theta({\bf z}) \cdot \partial_U\theta({\bf z}).
%\end{multline}	
%and if we set $\mathbf{z}_0$ to be a singular point of the theta divisor, we get:
\begin{equation}\label{eq:hirotadersingular}
{\partial_UH_{\bf z}}_{|\mathbf{z}=\mathbf{z}_0} = 2(\partial_U^2 \theta({\bf z_0}))(\partial^3_U \theta({\bf z_0}))
\end{equation}
 The quintic \eqref{eq:hirotadersingular} is a linear combination of the quintics $u_i \cdot F[\varepsilon]$, for $i=1,\dots,g$ and $\varepsilon \in \{0,1 \}^g$, so in principle we could try to proceed as for quadrics, and look for reducible quintics of the form $Q(U)\cdot T(U)$, where $\deg Q(U) = 2$ and $\deg T(U) =3$. However, this can be quite complicated, and it can be worth trying and computing directly a singular point of the theta divisor. Since this can be useful in general,  we discuss this in Section \ref{SingularPointsTheta}.
 
\subsection{Applications to the Schottky problem}\label{sec:schottky}
 
Up to now we have discussed the Torelli problem of reconstructing a smooth curve $C$ from a Riemann matrix $\tau$ of its Jacobian $J(C)$. Another fundamental question in this area is the Schottky problem \cite{Gru}, which asks, given a matrix $\tau \in \mathbb{H}_g$, whether this represents the Jacobian of a curve. This can be formulated in different ways with different possible solutions: see for example \cite{FGSM} for a very recent one. In particular, one of these was given by Krichever \cite{Kri} and Shiota \cite{Shi} via the KP equation. This solution can be formulated in terms of the Dubrovin threefold \cite[Section IV.4]{Dub81} by saying that $\tau\in \mathbb{H}_g$ represents a Jacobian if and only if the Dubrovin quartics \eqref{eq:dubrovinquartic} cut out a threefold. 

In particular, we can check that a matrix $\tau \in \mathbb{H}_g$ does not represent a Jacobian, by computing the quartics of Lemma \ref{lemma:quarticelimination} and then checking that they do not define a curve in $\mathbb{P}^{g-1}$. We verified this experimentally in Example \ref{example:schottky}.

\section{Numerical recovery}\label{dubrovin}

We can sum up the discussion of the previous section in the following algorithm. We have implemented it in Julia, which can be found at \url{https://turkuozlum.wixsite.com/tocj}. %and downloaded via \href{{https://wix.anyfileapp.net/dl?id=553246736447566b58312b65444764454f71746c6f4c526a78516851625868327639394d633635495a4a484f376d694e5a385462556e6f53794f494c31797441693461586a6d777233796b42425675323670456c32673d3d}}{\underline{this link}}. %{\color{red} insert a webpage link for the code.}
\vspace{5pt}

\begin{algorithm}[H]\label{alg:recovery}
	\KwIn{A  matrix $\tau\in \mathbb{H}_g$ representing the Jacobian of a non-hyperelliptic curve.}
	\KwOut{Quartics that cut out the canonical model of the algebraic curve $C$ whose Riemann matrix is $\tau$.}
	\KwSty{Step 1:} {Set up the linear system in~\eqref{eq:linearcondition} by computing the theta constants via the Julia package \texttt{Theta.jl}}.\\
	\KwSty{Step 2:} Solve the linear system in \eqref{eq:linearcondition}.\\
	\KwSty{Step 3:} Write the quartics~\eqref{eq:quarticelimination} and return them. \\
	\caption{Recovery through the Dubrovin threefold}
\end{algorithm}
\vspace{5pt}

The algorithm is straightforward and we can easily analyze its complexity in terms of the genus $g$. In Step 1, we need to evaluate $2^g\cdot (  {g(g+1)}/{2} + 1)$
theta constants, coming from the matrices $Q[\varepsilon]$ and the scalars $\hat{\theta}[\varepsilon]$. Then, in Step 2, we need to solve a $({g(g+1)}/{2} +1 ) \times  2^g$ linear system of maximal rank. Finally, in Step 3, we need to compute the quartics \eqref{eq:quarticelimination}, which involves the evaluation of 
$ 2^g \cdot {(g+3)(g+2)(g+1)g}/{24} $
theta constants. In our experiments, we considered examples, taken from the literature, up to genus $7$, so that the linear system of Step 2  is of relatively small size and it can be solved very quickly in Julia. What takes most of the time is the evaluation of the theta constants: the following table presents the approximated times to compute the theta constants in the examples below, with 12 digits of precision. In the table, $\partial^i$ indicates the order of the partial derivative of $\theta$ that we compute. The last column denotes the time needed to run the entire algorithm.

\begin{center}
\begin{tabular}{ | c | c | c | c | c |}
\hline
 genus & $\partial^0$ & $\partial^2$ & $\partial^4$ & total \\ \hline \hline
3 & 0.0009 sec & 0.001 sec & 0.002 sec & 5 sec \\ \hline
4 & 0.008 sec & 0.015 sec  & 0.02 sec & 11 sec \\ \hline  
5 & 0.07 sec & 0.15 sec  & 0.23 sec & 9 min \\ \hline  
6 & 2.1 sec & 4.2 sec  & 6.9 sec & 12 h \\ \hline  
7 & 6 sec & 8 sec  & 10 sec & 60 h \\ \hline  
\end{tabular}
\end{center}

\subsection{Computing the Singular Points}\label{SingularPointsTheta}
As we explained before, one of the advantages of the Dubrovin threefold is that it allows to recover the curve without computing a singular point of the theta divisor.   However, this is also a very useful method, if we manage to solve the transcendental system \eqref{eq:singularpointstheta}.
A Sage code that computes a singular point of the theta divisor in genus $4$ is presented in the article~\cite{ChuKumStu}. The basic idea is to solve system \eqref{eq:singularpointstheta} by numerical optimization, starting from a random input $z=a+\tau b$, where $a,b$ are real vectors with entries between $0$ and $1$.  In our implementation, we use the function \texttt{optimize.root} from the \texttt{SciPy} package. We call this function with the method \texttt{lm}, based on the Levenberg-Marquardt algorithm, which speeds up the computation substantially in comparison with the \texttt{hybr} method. The function \texttt{optimize.root} evaluates the partial derivatives~\eqref{eq:singularpointstheta} of the given function via estimating the limits of the function. Instead, we used the partial derivatives that is implemented in the Sage package \texttt{abelfunctions}~\cite{abelFunctions}, which gave more accurate results. In our experiments, it took about $30$ minutes for one singular point to be computed in the case of  genus $4$ and about $1.5$ hours in the case of genus $5$. 

\begin{remark}\label{rmk:changeofbasis}
	Before presenting our experiments, we observe that it is often convenient to work with an arbitrary basis of differentials instead of a normalized one as in \eqref{eq:adaptedbasis}. For such an arbitrary basis $\widetilde{\omega}_1,\dots,\widetilde{\omega}_g$, we consider the corresponding $g \times g$ period matrices
	\begin{equation}\label{eq:periodmatrices}
	\Pi_a \,=\, \left( \int_{a_j} \widetilde{\omega}_i \right)_{ij} \quad
	{\rm and} \qquad \Pi_b \,=\, \left( \int_{b_j} \widetilde{\omega}_i \right)_{ij}.
	\end{equation} 
	Then we obtain a normalized basis of differentials as in \eqref{eq:adaptedbasis} and the corresponding Riemann matrix by taking
	\begin{align}
	( \omega_1, \omega_2,\ldots,\omega_g)^T &= \Pi_a^{-1} \,
	( \widetilde{\omega}_1, \widetilde{\omega}_2,\ldots, \widetilde{\omega}_g)^T. \label{eq:basistransf} \\
	\tau &= \Pi_a^{-1}\Pi_b . \label{eq:riemannmatrix}
	\end{align}
\end{remark}

\subsection{Numerical experiments}

Finally, we present some examples illustrating our algorithm. 

In our experiments, we start with an explicit plane affine model for a non-hyperelliptic curve and possibily also its canonical model $C\subseteq \mathbb{P}^{g-1},$ and then we use the package \texttt{RiemannSurface} of Sage \cite{BruSijZot} to compute a Riemann matrix $\tau$ on which we run the Algorithm \ref{alg:recovery}. We then verify that the resulting quartics cut out the canonical curve we started with. We can do this explicitly in genus $3$, when the curve itself is a smooth plane quartic. In higher genera, we first verify that the quartics belong to the ideal of the curve by running the polynomial division algorithm, which return a remainder of zero, up to a certain numerical approximation. Furthermore, to verify that the quartics cut out the curve set-theoretically, we compute the intersection with an hyperplane in $\mathbb{P}^{g-1}$ by adding a random linear form and solving the resulting polynomial system via homotopy continuation. This is the primary computational method in numerical algebraic geometry, and we used the Julia implementation of  \texttt{HomotopyContinuation.jl}~\cite{BreTim}. This computation returns $2g-2$ solutions, confirming that the quartics cut out a curve of degree $2g-2$.

We also tried to recover the quadrics vanishing on the curve using the method of Section \ref{sec:quadrics}. We set up the problem of finding elements of the form $Q(U)^2$ in the space of quartics returned by Algorithm \ref{alg:recovery}, and we solved it again via \texttt{HomotopyContinuation.jl}. We could do this in genus $4$. %\textcolor{red}{THIS DID NOT WORK IN GENUS 5?}.

In genera $4$ and $5$, we could also compute singular points of the theta divisor, using the methods of Section \ref{SingularPointsTheta}. With these singular points, we could compute quadric and cubic equations for the curve, as described in Section \ref{sec:quadrics}.

%We used the
%{\tt SageMath} implementation due to Bruin et al.~\cite{BruSijZot} for computing the Riemann matrices and the period matrices. 

\begin{example}[Genus three]\label{example:Trott}
The Trott curve is a smooth plane quartic with affine model $ C = \{f(x,y) = 0\}$, where  \[ f(x,y ) = 12^2(x^4+y^4)-15^2(x^2+y^2)+350x^2y^2+81. \]
 In particular, this is already the canonical model, and the curve is of genus $3$ and not hyperelliptic. We compute a Riemann matrix using \texttt{RiemannSurface} in Sage \cite{BruSijZot}: in particular, the package uses the basis of differentials:
\[ \widetilde{\omega}_1 = \frac{1}{f_y}dx, \quad \widetilde{\omega}_2 = \frac{x}{f_y}dx, \quad  \widetilde{\omega}_3 = \frac{y}{f_y}dx, \] 
where $f_y$ denotes the derivative  $\frac{\partial f}{\partial y}$, and then it computes the period matrices $\Pi_a$ and $\Pi_b$ as in \eqref{eq:periodmatrices}. The corresponding normalized Riemann matrix $\tau$ is
\[
\begin{tiny}
\begin{pmatrix}
 1.06848368471179 + 0.723452867814272i & -0.305886633614305 + 0.123618182281837i & -0.160517941389541 - 0.206682546926085i \\
-0.305886633614305 + 0.123618182281837i &   0.776859918461210 + 1.25292663517205i & -0.626922516393387 - 0.289746911570334i \\
-0.160517941389541 - 0.206682546926085i & -0.626922516393387 - 0.289746911570334i &  0.376235735801471 + 0.484440302728207i \\
\end{pmatrix}.
\end{tiny}
\]
With this matrix, we set the linear system~\eqref{eq:linearcondition} by computing the theta constants appearing in the expressions of the system~\cite{julia}. As expected, this has an unique solution, up to a scalar multiplication, and we compute the corresponding quartic polynomial~\eqref{eq:quarticelimination}:
$$
\begin{scriptsize}
\begin{matrix}
(0.44055338231573327 - 0.11712521895532513i) u_1^4 + (2.094882287195226 + 7.879664904010854i) u_1^3 u_2
   \\- (5.316458517368645 - 1.4134300016965646i) u_1^3 u_3 + (61.49338091003442 - 16.348587918073555i) u_1^2 u_2^2 \\+ (27.505923029039046 + 105.6412469122926i) u_1^2 u_2 u_3 - (43.67750279381081 - 12.658628276584892i) u_1^2 u_3^2 \\
   - (0.20611709900405373 +0.7752863638524854i) u_1 u_2^3 + (142.137577271911 + 22.777083502115772i) u_1 u_2^2 u_3 \\
   + (101.16905240593528 + 146.6228999954985i) u_1 u_2 u_3^2 - (28.214458865117336 - 92.58798535078905i) u_1 u_3^3 \\
   - (0.06519271764094459 - 0.017332091034810038i) u_2^4 - (0.016856400506870983 + 0.8256030828721883i) u_2^3 u_3
   \\+ (64.66553470742735 + 38.49587006148285i) u_2^2 u_3^2 + (94.88897578016996 + 81.18194430047456i) u_2 u_3^3 \\
   + (33.080420780163195 + 41.521570514217885i) u_3^4.
\end{matrix}
\end{scriptsize}
$$
At a first glance, this might not look like the Trott curve. However, this equation is for the canonical model of $C$ with respect to a basis of normalized differentials $\omega_1,\omega_2,\omega_3$. If we go back to the differentials $\widetilde{\omega}_1,\widetilde{\omega}_2,\widetilde{\omega}_3$ via the change of coordinates in~\eqref{eq:basistransf}, we obtain the following quartic, after scaling the coefficients.  
%  First observation might be to intersect with a hypersurface as described in \ref{degree}. Using homotopy continuation, it can be observed that intersecting this quartic with the hypersurface $u_1+u_2+1$ which is a generic hypersurface will give $4$ solutions. Those $4$ solutions indicates that the degree of our curve is $4$. Also applying the change of basis described in .... will give the result:
$$
\begin{tiny}
\begin{matrix}
81u_1^4 + (1.2223597321441586\cdot10^{-13} - 9.454838005323456\cdot10^{-14}i)u_1^3u_2 \\ + (2.9124976279639876\cdot10^{-13} + 1.1282283371974781\cdot10^{-13}i)u_1^3u_3 - (225.00000000000017 - 4.607401108070593\cdot10^{-13}i)u_1^2u_2^2 \\ + (3.669767553538813\cdot10^{-13} - 3.017230893609506\cdot10^{-13}i)u_1^2u_2u_3 - (224.99999999999986 + 5.357443148919295\cdot10^{-13}i)u_1^2u_3^2 \\ - (4.1371303331328177\cdot10^{-13} - 3.5463573271644895\cdot10^{-13}i)u_1u_2^3 - (8.382029113614384\cdot10^{-13} + 3.97078497125283\cdot10^{-13}i)u_1u_2^2u_3 \\ + (7.725484571929981\cdot10^{-13} + 2.34428275709395\cdot10^{-13}i)u_1u_2u_3^2 - (8.239810406206657\cdot10^{-13} + 2.6625152861265\cdot10^{-13}i)u_1u_3^3 \\ + (143.99999999999918 - 7.607569271465399\cdot10^{-13}i)u_2^4 - (9.177234341211958\cdot10^{-13} - 8.304604428951876\cdot10^{-13}i)u_2^3u_3 \\ + (350.0000000000026 + 1.2750714694427922\cdot 10^{-12}i)u_2^2u_3^2 - (1.3400119435300388\cdot10^{-13} - 5.996042934502803\cdot10^{-13}i)u_2u_3^3 \\ + (143.99999999999895 - 5.357443148919295\cdot10^{-14}i)u_3^4.
\end{matrix}
\end{tiny}
$$
This is nothing but the quartic defining the Trott curve, up to an error of $10^{-12}$. In particular, we can recover the exact equation if we round up the coefficients to the nearest integer.  We emphasize that this example is treated slightly different than as it has been in \cite[Example 4.8]{AgoCelStu}.
\end{example}

Inspired by this example, we repeated the same experiments with 20 plane quartics with integer coefficients. The coefficients were bounded in absolute value by 100. We computed the period and the Riemann matrix with 53 bits of precision, we computed the theta constants with 12 digits of precision, and at the end we could recover the exact equation of the curve by rounding up the coefficients to the closest integer. Each experiment took approximately 4 seconds.

\begin{example}[Genus four] \label{ex:genus4}
Moving on to the case of genus 4, we consider the canonical curve
\begin{equation}\label{canonical4} 
C = \left\{ u_1 u_4 - u_2 u_3=0 \,,\,\, u_1^3 - u_2^3 - u_3^3 - u_4^3=0 \right\}.
\end{equation}
This has an affine plane model given by $\{ f(x,y) =0 \}$, where $f(x,y) = 1-x^3-y^3-x^3y^3$. We can recover the previous canonical model via the basis of differentials 
\[ \widetilde{\omega}_1 = -\frac{1}{f_y}dx,\quad \widetilde{\omega}_2=-\frac{x}{f_y}dx, \quad \widetilde{\omega}_3 = -\frac{y}{f_y}dx,\quad \widetilde{\omega}_4 = -\frac{xy}{f_y}dx. \] 

We can compute the $4 \times 4$ Riemann matrix $\tau$ via the plane model of the curve with the Sage package~\cite{BruSijZot}. This takes approximately 677 milliseconds for 53 bits, or about 16 digits, of precision.

%$
%\newline
%\tau =
%\begin{footnotesize}
%\left(
%\begin{matrix}
%    0.0816684226803939 + 1.25349451776915i &   - 0.760209337840853i \\
% - 0.760209337840853i &    -0.0816684226803936 + 1.25349451776915i \\
%    0.25 + 0.124170195056738i &     0.25 - 0.124170195056739i \\
%   -0.331668422680393 + 0.142753962855814i &    -0.25 - 0.636039142784115i \\
%\end{matrix} \right.
%\end{footnotesize}$
%$$\begin{footnotesize} \left.
%\begin{matrix}
%     0.25 + 0.124170195056738i &    -0.331668422680394 + 0.142753962855813i \\
%     0.25 - 0.124170195056738i &    -0.25 - 0.636039142784114i \\
%    0.760209337840853i &     0.24 - 0.636039142784115i \\
%    0.25 - 0.636039142784115i &      0.581668422680394 + 1.76536346549653i \\
%\end{matrix} \right)
%\end{footnotesize}$$
To reconstruct the canonical model of the curve back from $\tau$, we compute the $5$ quartics in~\eqref{eq:quarticelimination} by solving the linear system in~\eqref{eq:linearcondition}. %For a visual feast, we exhibit one of them: 

By Lemma \ref{lemma:quarticelimination}, these 5 quartics cut out the canonical curve \eqref{canonical4} after the basis change~\eqref{eq:basistransf}. We can first verify that the transformed quartics belong to the ideal of $C$ by the polynomial division algorithm. We did it in Sage, working over the complex field with 200 digits of precision. The coefficients of the remainder of the division algorithm were all of size $10^{-15}$.

Then, to verify that these quartic equations cut out the curve, we use the Julia package \texttt{HomotopyContinuation.jl}~\cite{BreTim}. We add a random linear form to the polynomial system of our 5 quartics and then \texttt{HomotopyContinuation.jl} returns $6$ solutions, which is what we expect from a curve of degree $6$ in $\mathbb{P}^3$. Moreover, again via homotopy continuation methods, we can find a quadratic polynomial $Q(U)$ such that $Q(U)^2$ is in the linear space generated by the 5 quartics, as in Section \ref{sec:quadrics}. %As seen below, this polynomial of degree $2$ does not seem like to our quadratic polynomial that is in the canonical ideal~\eqref{canonical4}.
After applying the change of basis in \eqref{eq:adaptedbasis} and rescaling we get the following expression for the quadric:%\textcolor{red}{CAN WE RESCALE THIS IN SUCH A WAY THAT THE COEFFICIENT OF $u_1u_4$ IS $1$?}:
$$
\begin{tiny}
\begin{matrix}
u_1   u_4 - (1.4829350744889013 \cdot 10^{-15} - 1.6682847904065378 \cdot 10^{-15}i)   u_1   u_2 \\- (3.5425309660018567 \cdot 10^{-15} + 6.403641669846521 \cdot 10^{-16}i)   u_1   u_3 (2.3679052278901118 \cdot 10^{-15} + 1.6728691462607347 \cdot 10^{-15}i)   u_1^2  \\ + (1.8423363133604865 \cdot 10^{-15} - 3.1265312370929112 \cdot 10^{-15}i)   u_2^2 - (1.0000000000000007 - 2.3672426468663847 \cdot 10^{-15}i)   u_2   u_3 \\- (1.739413822171687 \cdot 10^{-15} - 2.775912191360744 \cdot 10^{-15}i)  
u_2   u_4 + (3.3764916290020825 \cdot 10^{-15} - 6.001818720458345 \cdot 10^{-15}i)   u_3^2\\ - (3.777550457759975 \cdot 10^{-16} -1.4453231221486755 \cdot 10^{-15}i)   u_3   u_4 - (1.1268542006403671 \cdot 10^{-15} + 1.7990461567600933 \cdot 10^{-15}i)   u_4^2
\end{matrix}
\end{tiny}
$$

In genus four, we can also compute numerically a singular point of the theta divisor. We do it in Sage, as described in Section \ref{SingularPointsTheta} and we find the point:
\begin{equation*}
   {\bf z}_0=(0.75+0.54819629i, 0.75-0.54819629i,0.5+0.33618324i,0.75+0.2120130i).
\end{equation*}
The theta function and its derivatives vanish at this point up to 13 digits.
With this, we can compute the quadric $\partial^2_U\theta({\bf z}_0)$ and the cubic $\partial^3_U\theta({\bf z}_0)$. After the usual change of coordinates \eqref{eq:basistransf}, the quadric becomes 
%For this ${\bf z}_0$, all but one of the terms in \eqref{eq:hirotaquartics} vanish: $\,H_{{\bf z}_0} = 3(\partial^2_U\theta(\mathbf{z}_{0}))^2$. The quadratic polynomial $\partial^2_U\theta(\mathbf{z}_0)$ is the desired one as the Hirota quartics define the big Dubrovin threefold. We now apply the linear change of coordinate~\eqref{eq:adaptedbasis} to $\partial^2_U\theta(\mathbf{z}_0)$ as usual:
$$
\begin{tiny}
\begin{matrix}
 u_1u_4+(9.977112210552615\cdot 10^{-8} + 6.939529950175681\cdot 10^{-8}i)u_1^2 + (6.74409346713264\cdot 10^{-8} - 2.3021247380555947\cdot 10^{-15}
i)u_1u_2  \\
-(2.3274471968848503\cdot 10^{-8} + 5.037739720772648\cdot 10^{-8}i)u_1u_3+ (9.977111730319195\cdot 10^{-8} - 6.93953067335553\cdot 10^{-8}i)u_2^2  \\ 
-(0.9999999999999997 - 5.892639748496844\cdot 10^{-8}i)u_2u_3 + (2.3274465793326793\cdot 10^{-8} - 5.037739901039536\cdot 10^{-8}i)u_2u_4
\\
+ (3.9887820808350887\cdot 10^{-8} + 2.7942580923377634\cdot 10^{-8}i)u_3^2 + (1.3142269019133975\cdot 10^{-7} - 6.865574771844027\cdot 10^{-15}i)u_3u_4 
\\
+ (3.988781716962214\cdot 10^{-8} - 2.7942586271411155\cdot 10^{-8}i)u_4^2.
\end{matrix}
\end{tiny}
$$
which coincides with the quadric $u_1u_4-u_2u_3$ of \eqref{canonical4}, up to about 10 digits of precision. Instead, the cubic equation that we obtain is:
$$
\begin{tiny}
\begin{matrix}
u_1^3 - (1.0000000001244782 + 4.8426070737375934 \cdot 10^{-11}i)u_2^3 \\
-(1.0000000001244924 + 4.8417923378918607 \cdot 10^{-11}i)u_3^3- (1.0000000002161082 + 3.561043256396786 \cdot 10^{-11}i)u_4^3 \\
+ (4.5851667064228973 \cdot 10^{-11} + 3.883459414267029 \cdot 10^{-11}i)u_1^2u_2-(5.655054824682744 \cdot 10^{-11} - 2.0264603303086512 \cdot 10^{-11}i)u_1^2u_3\\
-(4.562873576948296 \cdot 10^{-11} + 1.0723355355290677 \cdot 10^{-10}i)u_1u_2^2+ (4.416105308034146 \cdot 10^{-11} + 6.485802277816666 \cdot 10^{-11}i)u_2^2u_4 \\
- (2.4095153783271284 \cdot 10^{-11} - 3.821774264257014 \cdot 10^{-11}i)u_2u_4^2- (2.1061545135963428 \cdot 10^{-11} + 3.9985852337772294 \cdot 10^{-11}i)u_3u_4^2 \\
+ (3.407561378730717 \cdot 10^{-11} - 7.066441338212726 \cdot 10^{-11}i)u_3^2u_4- (7.006134562951018 \cdot 10^{-11} - 9.313087455058934 \cdot 10^{-11}i)u_1u_3^2\\
+ (5.3419725864118215 - 2.3068027651869776i)u_1^2u_4- (5.341972586241411 - 2.3068027651197216i)u_1u_2u_3\\
- (0.861775203997493 - 1.4926384376919832i)u_1u_2u_4+ (0.8617752039804856 - 1.492638437827363i)u_2^2u_3 \\
-(0.8617752037076406 + 1.4926384378593538i)u_1u_3u_4 + (0.8617752038333869 + 1.4926384379123137i)u_2u_3^2 \\
- (2.396786904582313 - 2.79440847281551i)u_1u_4^2+ (2.3967869046516648 - 2.794408472851858i)u_2u_3u_4
\end{matrix}
\end{tiny}
$$
And we see that, with an approximation of 9 digits, this is the cubic $u_1^3-u_2^3-u_3^3-u_4^3$ of \eqref{canonical4}, plus with a linear combination of $u_i(u_1u_4-u_2u_3)$, for $i=1,2,3,4$.  \end{example}

\begin{example}\label{genus5example}
Let $C$ be the genus $5$ curve with an affine plane equation given by the polynomial $f(x,y)=x^2y^4 + x^4 + x + 3$. %The corresponding Riemann matrix $\tau$ can be calculated up using the \texttt{SageMath} package~\cite{BruSijZot}: \\
%
%$
%\tau=
%\begin{tiny}
%\left (
%\begin{matrix}
%  0.606845536499895 + 0.551363141414123i  & 0.129050724785112 - 0.301914886524400i \\
%  0.129050724785112 - 0.301914886524400i  & -0.234911158963853 + 1.19815149066113i  \\
% 0.115098558171014 + 0.0514312584207531i & -0.379202296886006 + 0.198016996468970i  \\
%0.492919888994140 + 0.00682403064655403i & -0.706761521228062 + 0.505080030041416i \\
%-0.232590742394293 - 0.256738092125412i &  0.101025144923035 - 0.0988712435720416i  
%\end{matrix} \right. \end{tiny}$ $$
%\begin{tiny}
%\left.
%\begin{matrix}
%0.115098558171014 + 0.0514312584207531i &0.492919888994140 + 0.00682403064655399i & -0.232590742394293 - 0.256738092125412i\\
% -0.379202296886007 + 0.198016996468970i & -0.706761521228061 + 0.505080030041416i & 0.101025144923035 - 0.0988712435720414i\\
% -0.393154463500105 + 0.551363141414122i & 0.492919888994141 + 0.00682403064655412i & -0.232590742394293 - 0.256738092125412i\\
% 0.492919888994140 + 0.00682403064655419i & -0.695539024056082 + 0.711658515659693i  & 0.516188714554510 - 0.214018280437827i\\
% -0.232590742394293 - 0.256738092125412i &  0.516188714554510 - 0.214018280437827i  & 0.341643728870296 + 0.881611263559767i
%\end{matrix}
%\right ).
%\end{tiny}
%$$
%Here we have 16 coming from the solutions of the linear system~\eqref{eq:linearcondition}. By using homotopy continuation we solve the zero dimensional polynomial system formed by these 16 quartics and $u_1+u_2+u_3+u_4+u_5$. There are 8 solutions to this system as expected. 
The differentials $\frac{1}{f_y}dx, \frac{x}{f_y}dx, \frac{xy}{f_y}dx, \frac{xy^2}{f_y}dx, \frac{x^2}{f_y}dx$ form a basis of the space of holomorphic differentials. The corresponding canonical model is given by the complete intersection of the three quadrics: 
\begin{equation*}
\label{canonicalGenus5}
u_4^2+u_5^2+u_2u_1+3u_1^2, \quad 
u_3^2-u_2u_4, \quad 
u_2^2-u_5u_1.
\end{equation*}
We compute the sixteen Dubrovin quartics~\eqref{eq:quarticelimination} and we check that, after the change of coordinates \eqref{eq:basistransf}, they belong to the ideal of $C$. We do this by polynomial division in Sage, over the complex field with 200 digits of precision as in Example~\ref{ex:genus4}. The coefficients of the remainder are of size $10^{-10}$.

We also check whether the quartics define a curve by adding a random linear form and solving the corresponding system via \texttt{HomotopyContinuation.jl}. We obtain $8$ solutions, which is nothing but the degree of our canonical genus 5 curve.  Here, one needs to increase the precision of about 15 digits while computing the Riemann matrix, which is required for computing the Dubrovin quartics. 

In this example, we could compute also singular points of the theta divisor as explained in Section \ref{SingularPointsTheta}. We computed three points $\mathbf{z}_1,\mathbf{z}_2,\mathbf{z}_2$ where the theta function and all its derivative vanish up to an error of $10^{-10}$. Then we obtain three quadrics $\partial^2_U\theta(\mathbf{z}_1),\partial^2_U\theta(\mathbf{z}_2),\partial^2_U\theta(\mathbf{z}_3)$, which, after the usual change of variables \eqref{eq:basistransf}, can be expressed as three independent linear combinations of the quadrics in \eqref{canonicalGenus5}, again up to an error of $10^{-10}$.
\end{example}

In the following examples, we push experimenting our methods to higher genera. 

\begin{example}[genus $6$ and $7$] Here we choose the curves of genus 6 and 7 as Wiman's sextic~\cite{Wim} and the butterfly curve~\cite{Fay}. Their respective plane affine equations are: 
\begin{align*}
     &x^{6}+y^{6}+1+(x^{2}+y^{2}+1)(x^{4}+y^{4}+1)=12x^{2}y^{2}, \\ 
     &\qquad \qquad \qquad \qquad \qquad x^{6}+y^{6}=x^{2}.
\end{align*}
We first compute their Riemann matrices numerically in Sage. Then, we estimate the corresponding 42 and 99 quartics~\eqref{eq:quarticelimination} in $\mathbb{P}^5$ and $\mathbb{P}^6$ respectively. Using the homotopy continuation method in Julia, we could verify that they define curves of degree 10 and 12 as expected. We point out that in these cases we needed to increase the precision in the Riemann matrix computation: 200 bits of precision in genus 6 and 500 bits in genus 7 were enough for the homotopy continuation computation to terminate. 
\end{example}

\begin{example}\label{example:schottky} 
Finally, we discuss some numerical experiments related to the Schottky problem, as in Section \ref{sec:schottky}. We choose 100 random Riemann matrices in genus 4, we computed the corresponding quartics as in Lemma \ref{lemma:quarticelimination}, we added a random linear form and we solved numerically the resulting system via \texttt{Homotopycontinuation.jl}. As expected, we found no solutions, confirming the fact that the quartics do not cut out a curve in $\mathbb{P}^{3}$. We expect that this circle of ideas would lead to an effective numerical solution to the Schottky problem, and we will investigate this in future work.
\end{example}

{\bf Acknowledgements}: We would like to thank Nils Bruins, Bernard Deconinck, Bernd Sturmfels and Andr\'{e} Uschmajew for their useful comments and their support.

\bigskip \medskip

\noindent
\footnotesize 
{\bf Authors' addresses:}

\smallskip

\noindent Daniele Agostini,  MPI-MiS Leipzig,
\hfill  {\tt daniele.agostini@mis.mpg.de}

\noindent 
 T\"urk\"u \"Ozl\"um \c{C}elik, Simon Fraser University,
\hfill {\tt turkuozlum@gmail.com}

%\noindent Bernard Deconinck,
%University of Washington,
%\hfill {\tt deconinc@uw.edu}

\noindent Demir Eken,
Bilkent University,
\hfill {\tt demir.eken@ug.bilkent.edu.tr}

\end{document}